\documentclass[12pt]{article}

\usepackage{epsfig}
\usepackage{amsfonts}
\usepackage{amssymb}
\usepackage{amsmath}   
\usepackage{amsthm}
\usepackage{latexsym}

\newtheorem{theorem}{Theorem}[section]
\newtheorem{lemma}[theorem]{Lemma}
\newtheorem{corollary}[theorem]{Corollary}

\newtheorem{definition}[theorem]{Definition}

\newtheorem{remark}[theorem]{Remark}

\newtheorem{system}{System}

\DeclareSymbolFont{AMSb}{U}{msb}{m}{n}
\DeclareMathSymbol{\N}{\mathbin}{AMSb}{"4E}
\DeclareMathSymbol{\Z}{\mathbin}{AMSb}{"5A}
\DeclareMathSymbol{\R}{\mathbin}{AMSb}{"52}
\DeclareMathSymbol{\Q}{\mathbin}{AMSb}{"51}
\DeclareMathSymbol{\I}{\mathbin}{AMSb}{"49}
\DeclareMathSymbol{\C}{\mathbin}{AMSb}{"43}

\title{Parameter estimation and asymptotic stability in stochastic filtering}

\author{Anastasia Papavasiliou\footnote{Current address: Department of Chemical Engineering, EQuad A217, Princeton, NJ 08544, USA.} \\ Columbia University}

\begin{document}
\maketitle

\begin{abstract}
In this paper, we study the problem of estimating a Markov chain $X$(signal) from its noisy partial information $Y$, when the transition probability kernel depends on some unknown parameters. Our goal is to compute the conditional distribution process ${\mathbb P}\{X_n|Y_n,\ldots,Y_1\}$, referred to hereafter as the {\it optimal filter}. Following a standard Bayesian technique, we treat the parameters as a non-dynamic component of the Markov chain. As a result, the new Markov chain is not going to be mixing, even if the original one is. We show that, under certain conditions, the optimal filters are still going to be asymptotically stable with respect to the initial conditions. Thus, by computing the optimal filter of the new system, we can estimate the signal adaptively. 
\end{abstract}

\bigskip
{\bf Key words:} nonlinear filtering, asymptotic stability, ergodic decomposition, Bayesian estimators.

\bigskip
{\bf AMS subject classifications:} 60G35, 93D20, 62F12.

\begin{section} {Introduction}
Stochastic filtering theory is concerned with the estimation of the distribution of a stochastic process at any time instant, given some partial information up to that time ({\it optimal filter}). The basic model usually consists of a Markov chain $X$ (also called the {\it state variable}) and possibly nonlinear observations $Y$ with observational noise $V$ independent of the signal $X$. In this case, the optimal filter is completely determined by the observations, the transition probability kernel, the distribution of the noise, and the initial distribution. One of the problems that often comes up in stochastic filtering is when one or more of these elements or, more generally, the model is not exactly known. In this paper, we study the case where the kernel depends on unknown parameters. 

The other important problem in stochastic filtering is how to compute the optimal filter. With the exception of very few cases (for example, linear Gaussian systems), an analytic solution does not exist and we have to resort to numerical methods. One of the most efficient schemes for the recursive computation of the optimal filter is the Interactive Particle Filter (or Sequential Monte Carlo), first suggested in \cite{Gordon} and \cite{Kitagawa}, independently. The idea is to approximate the optimal filter by an empirical distribution of particles which evolve in a way that imitates the evolution of the optimal filter. It has been shown that as the number of particles grows, the empirical distribution on these particles converges to the optimal filter, at every time instant (for theoretical results regarding the convergence of the Interacting Particle Filter see, for example, \cite{Del Moral-Guionnet}, \cite{Del Moral-Miclo}, \cite{Crisan-Del Moral-Lyons}, or \cite{Crisan-Doucet} for a comprehensive review).

These two problems are combined in the problem of adaptive estimation, i.e. how to estimate the parameters while computing the optimal filter. A natural idea is to treat the parameter as part of the state variable and then use some variation of the Interactive Particle Filter to compute the optimal filter (see \cite{West} for a historical perspective, as well as \cite{Storvik} and \cite{Doucet} for a more recent discussion). In this case, the Bayesian posterior distribution of the parameter is a marginal of the optimal filter. Even though there is plenty of numerical evidence showing that the posterior distribution of the parameter will converge to a delta function on the true value, this has not been proved yet, to the best of the writer's knowledge. The existing results on the consistency of estimators for the parameters of partially observed Markov chains concern other kinds of estimators (see \cite{Ryden} and references within for the case of Hidden Markov Models, i.e. partially observed Markov models with finite state space, or \cite{Moulines} for results regarding the consistency of the Maximum Likelihood Estimator for more general systems). We will show below that the Bayesian estimator is also consistent.

When we include the parameter in the state variable, the problem of adaptive estimation is clearly connected to the problem of asymptotic stability of the optimal filter with respect to its initial conditions. The true initial distribution on the second component of our state variable (the parameter) is of course a delta function on the true value, while we start with some prior distribution. It is also clear that the system is not mixing anymore, even if the signal $X$ is. The study of the asymptotic stability of the optimal filter is still an active area of research. In fact, many of the existing results for ergodic systems (\cite{Kunita 71},\cite{Stettner},\cite{Ocone-Pardoux}) have to be revised, since the recent discovery of a gap in a proof of \cite{Kunita 71} (see \cite{Lipster} and \cite{Budhiraja}). The question has been resolved for some cases as, for example, mixing systems (\cite{Atar-Zeitouni}) or particular cases of non-mixing systems (see \cite{Lipster},\cite{Chigansky}, \cite{Cerou} and \cite{Budhiraja-Ocone 99}). In this paper, we study the asymptotic stability of the optimal filters for systems that are not ergodic and whose ergodic components are actually mixing, under certain assumptions regarding the continuity of the kernel. 

Yet, as discussed above, usually we can only hope to compute an approximation to the optimal filter. This is always the case in adaptive estimation; even if the system is linear and Gaussian, the linearity is lost once we enter the parameters in the system. Since the error due to the unknown parameter disappears only asymptotically, we need a numerical scheme that converges uniformly with respect to time. The Interactive Particle Filters usually do not satisfy this condition if the system is not ergodic (see \cite{Del Moral-Guionnet}). A class of Particle Filters (the Monte-Carlo Particle Filter) that converges uniformly under relatively weak conditions is described in \cite{Del Moral}. These filters, though, are computationally much more expensive than the Interactive Particle Filters, since they require good sampling of the path space. In \cite{Papavasiliou}, we discuss a numerical scheme that converges uniformly to the optimal filter in question and is computationally relatively efficient.

The structure of the paper is the following: In section 2, we define the filtering problem we will be studying and state our main assumption, concerning the identifiability of the parameters. In section 3, we prove that the posterior distribution of the parameter will, indeed, converge to a delta function on the true value, under certain continuity conditions and for prior distributions on the parameter space whose mass in the neighborhood of the true value does not disappear ``too fast''.

In order to prove that the optimal filter is asymptotically stable, we first show that it is uniformly continuous with respect to the parameters. The uniform continuity of the optimal filter has also been studied in \cite{Oudjane}, under the assumption that the kernels are mixing. In that case, it is shown that if the parameter is fixed to a value different than the true one, the total error in the optimal filter will be uniformly bounded by the supremum of the step errors -- i.e. the errors made when the parameters are different just at the last time step -- multiplied by some constant that depends on the mixing constant. Thus, the uniform continuity of the step errors implies the uniform continuity of the optimal filter. In section 4, we review these results and give some sufficient conditions for the uniform continuity of the step error to hold, that are relatively easy to check. Finally, in section 5, we prove the main result: the asymptotic stability of the optimal filter, with respect to the initial conditions, for systems that come up when we include the parameters of an otherwise mixing system, in the state variable.
 
\end{section}

\section{Definitions and Assumptions}

Let $E$ be a Polish space, i.e., a complete separable metric space and let us denote by ${\mathcal B}(E)$ its Borel $\sigma$-field. We study the asymptotic behavior of the conditional distribution of a Markov chain $\{X_n\}$ taking values in $E$ given some noisy partial information, when the kernel depends on an unknown parameter $\theta$. More specifically, we study the optimal filter of the following system, which we will refer to as:

\bigskip
\begin{system}
\label{system A}
Let $\{X_n\}$ be a homogeneous Markov chain taking values in $(E,{\mathcal B}(E))$. Let $\mu$ be its initial distribution and $K_\theta$ its transition probability kernel depending on a parameter $\theta\in\Theta$. Furthermore, we assume that for each $\theta\in\Theta$, $K_\theta$ is Feller and mixing, i.e. there exists a constant $0<\epsilon_\theta\leq 1$ and a nonnegative measure $\lambda_\theta\in{\mathcal{M}}^{+}(E)$ (${\mathcal{M}}^{+}(E)$ being the set of finite nonnegative measures on $E$), such that 
\begin{equation}
\label{mixing condition}
\epsilon_\theta\lambda_\theta(A) \leq K_\theta(x,A) \leq \frac{1}{\epsilon_\theta}\lambda_\theta(A),\ \ \ \forall x\in E \ and\ \forall A\in{\mathcal B}(E).
\end{equation}
The observation process is defined by 
\begin{displaymath}
Y_n=h(X_n) + V_n,
\end{displaymath} 
where $V_n$ are i.i.d. Gaussian(0,$\sigma^2$) $\R^p$-valued random variables independent of $X$ and ${h:E\rightarrow\R^p}$ is a bounded continuous function. We denote by $g$ the Gaussian probability density function of the observational noise.
\end{system}
\bigskip

In practice, the parameter space $\Theta$ is usually Euclidean. More generally, we assume that it is a Polish space, with metric $d_\Theta(\cdot,\cdot)$. Most problems are given in the form of System 1. Following a standard Bayesian technique, we rewrite the system, so that the parameter becomes part of the Markov chain, whose transition probability kernel is now completely known:

\bigskip
\begin{system}
\label{system B}
Suppose that $\{\tilde{X}_n=(X_n,\theta_n)\}$ is an $E\times\Theta$-valued homogeneous Markov chain, with transition probability 
\begin{displaymath}
\tilde{K}((x^\prime,\theta^\prime),dx\otimes d\theta)=K_\theta(x^\prime,dx)\otimes\delta_{\theta^\prime}(d\theta).
\end{displaymath}
and initial distribution $\mu\otimes u$, where $K_\theta$ is Feller and mixing in the sense of (\ref{mixing condition}). The observation process is defined by
\begin{displaymath}
Y_n = \tilde{h}(\tilde{X}_n) + V_n,
\end{displaymath}
where $\tilde{h}(\tilde{x})=h(x)$ and $\tilde{x}=(x,\theta)$ and the process $\{V_n\}_{n>0}$ is defined as above.
\end{system}
\bigskip

System 2 can be thought of as a generalization of System 1 and will be the main object of study in this paper. Our goal is to show the asymptotic stability of the optimal filter of this system with respect to its initial distribution.

The canonical space of the Markov chain $X$ with kernel $K_\theta$ and initial distribution $\mu$ is denoted by  $(\Omega_1 = E^{\mathbb N},({\mathcal F}^{(X)}_n)_{n\geq 0}, P_{\mu,\theta})$, where ${\mathcal F}^{(X)}_n = \sigma(X_0, X_1,\dots,X_n)$ is the $\sigma$-algebra constructed by the random variables $X_0,X_1,\dots,X_n$. Similarly, the observation process is defined on the canonical space $(\Omega_2 = ({\mathbb R}^p)^{\mathbb N},({\mathcal F}^{(Y)}_n)_{n\geq 0}, Q_{\mu,\theta})$, where ${\mathcal F}^{(Y)}_n = \sigma(Y_1,\dots,Y_n)$. The law of the observation process $Q_{\mu,\theta}$ is given by 
\begin{displaymath}
Q_{\mu,\theta}(dy_{k_1},\dots,dy_{k_n}) = \int_{E^{\otimes n}} \prod_{i=1}^n g(y_{k_i} - h(x_{k_i}))P_{\mu,\theta}(dx_{k_1},\dots,dx_{k_n})dy_{k_1}\dots dy_{k_n},
\end{displaymath}
for any $n\geq 0$ and $k_1,\dots,k_n\in\R_{+}$, where $E^{\otimes n} = E\times \dots \times E$ is the product space of n copies of $E$. Also, by $Q^n_{\mu,\theta}$ we denote the restriction of the measure $Q_{\mu,\theta}$ to the sigma algebra ${\mathcal F}^{(Y)}_n$.

We can now define the pair process $(X,Y)$ on the space $({\Omega} = \Omega_1\times\Omega_2, ({\mathcal F}_n = {\mathcal F}^{(X)}_n \times{\mathcal F}^{(Y)}_n)_{n\geq 0}, {\mathbb P}_{\mu,\theta})$, where the measure ${\mathbb P}_{\mu,\theta}$ is such that its marginal distributions with respect to $X$ and $Y$ are $P_{\mu,\theta}$ and $Q_{\mu,\theta}$ respectively. It is not hard to show that this measure exists (see, for example, \cite{Del Moral}). We will denote the expectation with respect to ${\mathbb P}_{\mu,\theta}$ by ${\mathbb E}_{\mu,\theta}$.

Similarly, we define the triplet $(X,Y,\theta)$ on the space $(\tilde{\Omega} = {\Omega}\times\Omega_3, (\tilde{\mathcal F}_n = {\mathcal F}_n \times\sigma(\theta))_{n\geq 0}, \tilde{\mathbb P}_{\mu,u})$, where $\theta$ is a $\Theta$-valued random variable defined on $(\Omega_3, \sigma(\theta), u)$ and the marginals of $\tilde{\mathbb P}_{\mu,u}$ on $(X,Y)$ and $\theta$ respectively are $\int_\Theta {\mathbb P}_{\mu,\theta} u(d\theta)$ and $u$. We will denote the expectation with respect to $\tilde{\mathbb P}_{\mu,u}$ by $\tilde{\mathbb E}_{\mu,u}$.

Furthermore, we denote by $\Psi^\theta_n(\mu)$ and $\Phi_n(\mu\otimes u)$ the optimal filters for Systems \ref{system A} and \ref{system B}, with initial distributions $\mu$ and $\mu\otimes u$, respectively, defined as the posterior distribution of the state variable given the observations. The name ``optimal filters'' is due to the fact that they are the best estimators adapted to the available information (the $\sigma$-algebra constructed by the observations), with respect to the $L_2$-norm. They are random measures on the space $E$ and $E\times \Theta$ respectively, defined as follows: for every $f\in{\mathcal C}_b(E\times\Theta)$,

\begin{eqnarray}
\label{Phi_n}
\nonumber\Phi_n(\mu\otimes u)(f) &=& \tilde{\mathbb E}_{\mu,u}[f(X_n,\theta)|Y_n,\ldots,Y_1] = \\ \nonumber\\
&=& \frac{\int_\Theta\int_{E^{\otimes n}} f(x_n,\theta)\prod^{n}_{k=1} g(Y_k-h(x_k))P_{\mu,\theta}(dx_1,\dots dx_n)u(d\theta)} {\int_\Theta\int_{E^{\otimes n}}\prod_{k=1}^n g(Y_k-h(x_k)) P_{\mu,\theta}(dx_1,\dots,dx_n)u(d\theta)}.
\end{eqnarray}

\noindent
Similarly, $\forall f^\prime\in{\mathcal C}_b(E)$, 

\begin{eqnarray}
\label{Psi_n}
\nonumber \Psi^\theta_n(\mu)(f^\prime) &=& {\mathbb E}_{\mu,\theta}[f^\prime(X_n)|Y_n,\ldots,Y_1] = \\ \nonumber\\
&=& \frac{\int_{E^{\otimes n}} f^\prime(x_n)\prod_{k=1}^n g(Y_k-h(x_k))P_{\mu,\theta}(dx_1,\dots,dx_n)}{\int_{E^{\otimes n}}\prod_{k=1}^n g(Y_k-h(x_k)) P_{\mu,\theta}(dx_1,\dots,dx_n)}.
\end{eqnarray}

Clearly, $\Psi_n^\alpha(\mu)$ is the marginal of $\Phi_n(\mu\otimes\delta_\alpha)$ with respect to the state variable $X$, since $\Phi_n(\mu\otimes\delta_\alpha)(f)=\Psi^\alpha_n(\mu)(f_\alpha)$, where we used the notation $f_\alpha(x) = f(x,\alpha)$. Similarly, we define $\Psi^u_n(\mu)$ as the marginal of $\Phi_n(\mu\otimes u)$ with respect to $X$, i.e., 
\[\Psi^u_n(\mu)(dx) = \int_\Theta \Phi_n(\mu\otimes u)(dx,d\theta).\] 
We want to find sufficient conditions for $\Psi^u_n(\mu)$ to be asymptotically stable with respect to the initial distribution $u$. First, we need to check the ``identifiability'' of the parameters.

The mixing condition (\ref{mixing condition}) implies the ergodicity of the signal. We denote by $\mu_\theta$ the limiting distribution of the Markov chain whose transition kernel is $K_\theta$. The measure $\mu_\theta$ is uniquely defined in this way, as a result of the ergodic property of the signal. By $\nu_\theta$ we denote the limiting distribution of the observation process corresponding to parameter value $\theta$. It is easy to check that this is also uniquely defined by: 
\begin{equation}
\label{nu}
\nu_\theta  = (\mu_\theta\circ h^{-1})\ast g.
\end{equation}

We define an equivalence relation on the parameter space as follows:
\begin{equation}
\label{equivalence}
\alpha\sim\beta\Leftrightarrow\mu_\alpha\circ h^{-1} = \mu_\beta\circ h^{-1}
\end{equation}
Since $g$ is the Gaussian probability density function, the following definition is equivalent to (\ref{equivalence}):
\begin{equation}
\label{equivalence2}
\alpha\sim\beta\Leftrightarrow\nu_\alpha = \nu_\beta.
\end{equation}

We assume that there is no pair of equivalent points in the parameter space. This implies that the observation processes corresponding to two different parameter values are mutually singular. Otherwise, we might not be able to tell two parameter values apart by looking at the observations. A trivial example is when $h$ is constant. Problems can also arise when $h$ is symmetric. 

\bigskip
\noindent
{\bf Identifiability condition:}{\it\ By saying that ``the identifiability condition holds'' or that ``the parameters are identifiable'', we mean that $\alpha\neq\beta$ implies $\alpha\not\sim\beta$.}

\begin{remark}
The assumption that the noise $\{V_n\}_{n>0}$ in the observation process is Gaussian is not really necessary. It just implies the equivalence between definitions (\ref{equivalence}) and (\ref{equivalence2}). For arbitrary observational noise, one just has to define the equivalence relation using (\ref{equivalence2}) and all the following results will also hold, with the appropriate changes in notation.
\end{remark}

\begin{section}{Consistency of Bayesian estimator}

In this section, we study the behavior of the posterior distribution of the parameter, given the observations. We show that under certain conditions and for almost every realization $y=\{y_1,\dots,y_n,\dots\}$ of the observation process described by System \ref{system A} corresponding to a fixed value $\alpha$ of the parameter ($\theta=\alpha$), the posterior distribution of the parameter $\tilde{\mathbb P}_{\mu,u}(\theta | y)$, where $u$ is the prior distribution, is a delta function on $\alpha$. 

If we only assume that the identifiability condition holds, we can show that for $u$-a.e. $\alpha$,
\begin{equation}
\label{posterior eq}
Q_\alpha\{\tilde{\mathbb P}_{\mu,u}(\theta|y) = \delta_\alpha\} = 1,
\end{equation}
meaning that for almost every realization $\alpha$ of the parameter value with respect to the probability measure $u$  and almost every realization $y$ of the observation process with respect to the measure $Q_\alpha$, the posterior distribution of the parameters given the observation is going to be a delta function on $\alpha$. Equivalently, this can be stated as follows:
\begin{equation*}
u{\big\{}\alpha\in\Theta: Q_\alpha\{y: \tilde{\mathbb P}_{\mu,u}(\theta|y)=\delta_\alpha\}=1{\big\}}=1.
\end{equation*}
This is a consequence of a lemma by Gl\"{o}tzl and Wakolbinger \cite{Glotzl}, where they use the notion of ergodic decomposability. It is a promising result, but we actually need something stronger. Since we assume that there exists a ``true value'' for the parameter which is fixed but unknown, we have to be sure that whatever this value is, the result will hold. So, we would like prove (\ref{posterior eq}) for every possible value of the parameter. 

Below, we prove a result for the posterior distribution that is slightly weaker than (\ref{posterior eq}), for any ``reasonable prior'' (this will be made more precise below) and under some additional assumptions. This result, however, will be sufficient for our purposes.

\begin{lemma}
\label{posterior lemma 1}
Let $Y$ be as described in System \ref{system A}. We further assume that
\begin{itemize}
	\item{ the prior distribution $u$ on the parameter space is such that there exist a sequence $\epsilon_n\downarrow0$ and a function $p:{\mathbb N}\rightarrow[1,\infty)$ with the following properties:
\begin{equation}
\label{p(n)}
\frac{p(n)}{n}\rightarrow0,\ {\rm as}\ n\rightarrow\infty
\end{equation}
and
\begin{equation}
\label{condition on prior}
\sup_{n\uparrow\infty}{\mathbb E}_{\mu,\alpha}[(\frac{\sup_{\theta\in{\mathcal N}_{\epsilon_n}(\alpha)}\frac{dQ^n_{\mu,\alpha}}{dQ^n_{\mu^\prime,\theta}}(Y_n,\dots,Y_1)}{u({\mathcal N}_{\epsilon_n}(\alpha))})^{\frac{1}{p(n)}}]<+\infty.
\end{equation}}
	\item{For each $\eta>0$, there exists an $\epsilon>0$ and an $I_\eta>0$ such that 
\begin{equation}
\label{full LDE}
\limsup_{n\rightarrow\infty}\frac{1}{n}\log\int_{{\mathcal N}_\eta(\alpha)^c}{\mathbb P}_{\mu^\prime,\theta}(L_n(Y)\in{\mathcal B}_\epsilon(\nu_\alpha))u(d\theta)\leq-I_\eta<0,
\end{equation}
where $L_n(Y)$ is the empirical measure of the observation process up to time $n$, i.e. $L_n(Y) = \frac{1}{n}\sum_{k=1}^n \delta_{Y_k}$.}
\end{itemize}
Then, for every $\eta>0$,
\begin{equation}
\label{posterior}
\lim_{n\rightarrow\infty}{\mathbb E}_{\mu,\alpha}\tilde{\mathbb P}_{\mu^\prime,u}(\theta\in {\mathcal N}_\eta(\alpha)^{c} |Y_n,\dots,Y_1) =0.
\end{equation}
Note that we have used the notation ${\mathcal N}_\rho(\alpha)$ to denote the ball of radius $\rho$ and center $\alpha$ with respect to the metric $d_\Theta$, for any $\rho>0$. Similarly,  by ${\mathcal B}_\epsilon(\nu_\alpha)$ we denote the ball of radius $\epsilon$ and center $\nu_\alpha$ with respect to the L\'evy-Prohorov metric. Also, by $A^{c}$ we denote the complement of any set $A$.
\end{lemma}

The condition on the prior says that the mass around the true value $\alpha$ should not go to zero ``too quickly'' and how quickly is that will depend on how fast the measure $Q_\theta$ approaches $Q_\alpha$ when $\theta$ goes to $\alpha$. This condition becomes more clear when applied to specific models.

The proof of lemma \ref{posterior lemma 1} is given below:

\begin{proof}[Proof of lemma \ref{posterior lemma 1}]
Let us fix an $\eta>0$. Then, we can choose $\epsilon>0$ so that (\ref{full LDE}) holds.
We break the expectation in two parts: 
\begin{equation*}
{\mathbb E}_{\mu,\alpha}\tilde{\mathbb P}_{\mu^\prime,u}(\theta\in {\mathcal N}_\eta(\alpha)^{c} |Y_n,\dots,Y_1) = A_n + B_n,
\end{equation*}
where
\[ A_n = {\mathbb E}_{\mu,\alpha}[{\mathbf 1}_{{\mathcal B}_\epsilon(\nu_\alpha)}(L_n(Y))\tilde{\mathbb P}_{\mu^\prime,u}(\theta\in {\mathcal N}_\eta(\alpha)^{c} |Y_n,\dots,Y_1)]\] 
and
\[ B_n = {\mathbb E}_{\mu,\alpha}[{\mathbf 1}_{({\mathcal B}_\epsilon(\nu_\alpha))^c}(L_n(Y))\tilde{\mathbb P}_{\mu^\prime,u}(\theta\in {\mathcal N}_\eta(\alpha)^{c} |Y_n,\dots,Y_1)].\]
Clearly
\[ \lim_{n\rightarrow\infty}B_n \leq \lim_{n\rightarrow\infty}{\mathbb P}_{\mu,\alpha}(L_n(Y)\in({\mathcal B}_\epsilon(\nu_\alpha))^c)=0,\ \forall\epsilon>0 \]
by Birkhoff's ergodic theorem. It remains to show that $\lim_{n\rightarrow\infty}A_n = 0$. We write
\begin{eqnarray*}
A_n &=& \tilde{\mathbb E}_{\mu^\prime,u}[{\mathbf 1}_{({\mathcal B}_\epsilon(\nu_\alpha))^c}(L_n(Y))\cdot\tilde{\mathbb E}_{\mu^\prime,u}({\mathbf 1}_{({\mathcal N}_\eta(\alpha))^c}(\theta)|Y_n,\dots,Y_1)\frac{dQ^n_{\mu,\alpha}}{dQ^n_{\mu^\prime,u}}(Y_n,\dots,Y_1)] = \\
&=& \tilde{\mathbb E}_{\mu^\prime,u}[{\mathbf 1}_{({\mathcal B}_\epsilon(\nu_\alpha))^c}(L_n(Y))\cdot{\mathbf 1}_{({\mathcal N}_\eta(\alpha))^c}(\theta)\frac{dQ^n_{\mu,\alpha}}{dQ^n_{\mu^\prime,u}}(Y_n,\dots,Y_1)] \leq \\
&\leq& \tilde{\mathbb E}_{\mu^\prime,u}[({\mathbf 1}_{({\mathcal B}_\epsilon(\nu_\alpha))^c}(L_n(Y)){\mathbf 1}_{({\mathcal N}_\eta(\alpha))^c}(\theta))^{p(n)+1}]^{\frac{1}{p(n)+1}}\tilde{\mathbb E}_{\mu^\prime,u}[(\frac{dQ^n_{\mu,\alpha}}{dQ^n_{\mu^\prime,u}})^{\frac{p(n)+1}{p(n)}}]^{\frac{p(n)}{p(n)+1}} = \\
&=& (\tilde{\mathbb E}_{\mu^\prime,u}[{\mathbf 1}_{({\mathcal B}_\epsilon(\nu_\alpha))^c}(L_n(Y)){\mathbf 1}_{({\mathcal N}_\eta(\alpha))^c}(\theta)])^{\frac{1}{p(n)+1}}({\mathbb E}_{\mu,\alpha}[(\frac{dQ^n_{\mu,\alpha}}{dQ^n_{\mu^\prime,u}})^{\frac{1}{p(n)}}])^{\frac{p(n)}{p(n)+1}}.
\end{eqnarray*}
To go from line 1 to line 2 above we used the tower property of the expectations and the adaptiveness of the terms outside the conditional expectation to the $\sigma$-algebra $\sigma(Y_n,\dots,Y_1)$. To go from line 2 to line 3 we applied the Holder inequality with $p=p(n)+1$ and $q=\frac{p(n)+1}{p(n)}$.
The first term can be written as 
\[ \exp(\frac{n}{p(n)+1}\frac{1}{n}\log\int_{({\mathcal N}_\eta(\alpha))^c}{\mathbb P}_{\mu^\prime,\theta}(L_n(Y)\in{\mathcal B}_\epsilon(\nu_\alpha))u(d\theta)) \]
and it goes to zero as $n\rightarrow\infty$, by (\ref{p(n)}) and (\ref{full LDE}). The second term will be uniformly bounded, since
\begin{eqnarray*}
{\mathbb E}_{\mu,\alpha}[(\frac{dQ^n_{\mu,\alpha}}{dQ^n_{\mu^\prime,u}})^{\frac{1}{p(n)}}] &=& {\mathbb E}_{\mu,\alpha}[(\frac{dQ^n_{\mu,\alpha}}{\int_\Theta dQ^n_{\mu^\prime,\theta} u(d\theta)})^{\frac{1}{p(n)}}] \leq \\
{\mathbb E}_{\mu,\alpha}[(\frac{dQ^n_{\mu,\alpha}}{\int_{{\mathcal N}_{\epsilon_n}(\alpha)} dQ^n_{\mu^\prime,\theta} u(d\theta)})^{\frac{1}{p(n)}}] &\leq& {\mathbb E}_{\mu,\alpha}[(\frac{\sup_{\theta\in{\mathcal N}_{\epsilon_n}(\alpha)}\frac{dQ^n_{\mu,\alpha}}{dQ^n_{\mu^\prime,\theta}}(Y_n,\dots,Y_1)}{u({\mathcal N}_{\epsilon_n}(\alpha))})^{\frac{1}{p(n)}}].
\end{eqnarray*}
Thus, $\lim_{n\rightarrow\infty}A_n = 0$, which completes the proof.
\end{proof}

From the proof of lemma \ref{posterior lemma 1}, one can easily get an estimate for the rate of convergence of (\ref{posterior}). If the Markov chain for $\theta=\alpha$ satisfies the LDP, then
\begin{equation}
\label{posterior rate}
\limsup_{n\rightarrow\infty}\frac{p(n)+1}{n}\log {\mathbb E}_{\mu,\alpha}\tilde{\mathbb P}_{\mu^\prime,u}(\theta\in {\mathcal N}_\eta(\alpha)^{c} |Y_n,\dots,Y_1) \leq -I_\eta <0.
\end{equation}

Note that if the prior $u$ puts positive mass on the true value $\alpha$, i.e. $u(\{\alpha\})>0$, we can choose $\epsilon_n\equiv0$ and $p(n)=n$. Then, the rate of convergence will be exponential (in the sense of (\ref{posterior rate})).

We would now like to find conditions for (\ref{full LDE}) to hold that only depend on the properties of the kernels $K_\theta$. If $u({\mathcal N}_\eta(\alpha)^c)=0$ it holds trivially, so we assume that $u({\mathcal N}_\eta(\alpha)^c)>0$. Let us also assume that for each $\theta\in\Theta$, the Markov chain satisfies the LDP with rate function $I_\theta$. This is really just a property of the kernels $K_\theta$ (see, for example, \cite{Dupuis-Ellis}, for the properties the kernel has to satisfy in order for the Markov chain to satisfy the LDP). Then, by the contraction principle, the observation process $\{Y_n\}_{n>0}$ will also satisfy the LDP with a good rate function $J_\theta$, given by
\begin{eqnarray}
\label{J theta}
&J_\theta(\nu):= \inf\{\inf_{\{K\in{\mathcal T}: \mu K=\mu\}} \int_E\int_E \log{\frac{dK(x,\cdot)}{dK_\theta(x,\cdot)}(y)}K(x,dy)\mu(dx) + \\
\nonumber & + \inf_{\{v^\prime\in{\mathcal P}(\R^p)\}} \int_{\R^p} \log{\frac{dv^\prime}{dv}(x)}v^\prime(dx)
;\ \mu\in{\mathcal P}(E), v\in{\mathcal P}(\R^p)\ {\rm and}\ \nu = (\mu\circ h^{-1})\ast v\}
\end{eqnarray}
where ${\mathcal T}$ is the class of all transition kernels on $E$.
To show (\ref{full LDE}), we have to show that
\[ \limsup_{n\rightarrow\infty}\sup_{\theta\in{\mathcal N}_\eta(\alpha)^c}\frac{1}{n}\log{\mathbb P}_{\mu^\prime,\theta}(L_n(Y)\in{\mathcal B}_\epsilon(\nu_\alpha))\leq -I_\eta <0. \]
If the parameter space is compact and the distribution ${\mathbb P}_{\mu^\prime,\theta}$ is continuous with respect to $\theta$, we can interchange the limit and the supremum and, consequently, it suffices to show that for each $\theta\in {\mathcal N}_\eta(\alpha)^c$,
\[ \limsup_{n\rightarrow\infty}\frac{1}{n}\log{\mathbb P}_{\mu^\prime,\theta}(L_n(Y)\in{\mathcal B}_\epsilon(\nu_\alpha))\leq -I_\eta <0. \]
Since the observation process satisfies the LDP for each $\theta$, we know that
\begin{equation*}
\limsup_{n\rightarrow\infty}\frac{1}{n}\log{\mathbb P}_{\mu^\prime,\theta}(L_n(Y)\in{\mathcal B}_\epsilon(\nu_\alpha))\leq -J_\theta(\overline{{\mathcal B}_\epsilon(\nu_\alpha)}),\ \forall\theta\in\Theta, \end{equation*}
where by $\overline{{\mathcal B}_\epsilon(\nu_\alpha)}$ we denote the closure of ${\mathcal B}_\epsilon(\nu_\alpha)$. We need to find an $\epsilon>0$, such that 
\begin{equation}
\label{J}
\inf_{\theta\in{\mathcal N}_\eta(\alpha)^c} J_\theta(\overline{{\mathcal B}_\epsilon(\nu_\alpha)}) >0.
\end{equation}
Let us also assume that the rate function $J_\theta$ is continuous with respect to $\theta$. Then, the compactness of the parameter space and the properties of rate functions imply that (\ref{J}) will be true if for each $\theta\in{\mathcal N}_\eta(\alpha)^c$, $\nu_\theta\notin\overline{{\mathcal B}_\epsilon(\nu_\alpha)}$. Equivalently, we ask that 
\begin{equation}
\label{open mapping}
\forall\eta>0,\ \exists\epsilon>0: (L(\nu_\theta,\nu_\alpha)\leq\epsilon \Rightarrow d_\Theta(\theta,\alpha)< \eta),
\end{equation}
where by $L(\cdot,\cdot)$ we denote the L\'evy-Prohorov metric. This holds if the mapping from the parameter space $\Theta$ to the space of limiting distributions of the observation process is open (i.e. open sets are matched to open sets).

Note that the continuity of the kernel $K_\theta$ with respect to the parameter $\theta$ implies the continuity of ${\mathbb P}_{\mu^\prime,\theta}$ with respect to $\theta$. So, we have found conditions for (\ref{full LDE}) to hold that only involve properties of the kernels. These are summarized in the following 

\begin{lemma}
\label{posterior lemma 2}
Suppose that the parameter space $\Theta$ is compact and the following hold:
\begin{itemize}
	\item{for each $\theta\in\Theta$, the observation process $\{Y_n\}_{n>0}$ satisfies the LDP with rate function $J_\theta$. }
	\item{The mapping $\Theta\ni\theta\mapsto{\nu_\theta}\in{\mathcal P}({\mathbb R}^p)$, is an open mapping, in the sense of (\ref{open mapping}).}
	\item{The kernel $K_\theta$ and the rate function $J_\theta$ are continuous with respect to $\theta$.}
\end{itemize}
Then, (\ref{full LDE}) will also hold.
\end{lemma}

\end{section}

\section{Uniform continuity of the optimal filter}

In this section, we study uniform continuity of the optimal filter, with respect to the parameter. This problem was first studied in \cite{Oudjane}, where it was shown that if the kernel is mixing and the step errors converge to zero uniformly, then the total error is also going to converge to zero uniformly. We review these results here and give some sufficient conditions for the step error to be uniformly continuous. We first review the definition of the ``Hilbert projective metric'' in the space of measures:

\begin{definition}
\label{Hilbert metric}
Two nonnegative measures $\mu,\mu^\prime \in{\mathcal M}^{+}(E)$ are {\bf comparable}, if $\alpha \mu \leq \mu^\prime \leq \beta \mu$ for suitable positive scalars $\alpha$ and $\beta$. The {\bf Hilbert projective metric} on ${\mathcal M}^{+}(E)$ is defined as
\begin{equation}
h(\mu,\mu^\prime) := 
\left\{ \begin{array}{c}
	\log\frac{\sup_{A: \mu^\prime(A)>0} \frac{\mu(A)}{\mu^\prime(A)}}{\inf_{A: \mu^\prime(A)>0} \frac{\mu(A)}{\mu^\prime(A)}} = \log(\|\frac{d\mu}{d\mu^\prime}\|_{\infty}\|\frac{d\mu^\prime}{d\mu}\|_{\infty}), \\  \\ 0, \\ \\
	+\infty,
\end{array}
\begin{array}{c}
if\ \mu\ and\ \mu^\prime\ are\ comparable \\ \\
if\ \mu=\mu^\prime\equiv 0 \\ \\
otherwise 
\end{array}\right.
\end{equation}

The kernel norm corresponding to the Hilbert metric is called the {\bf Birkhoff contraction coefficient} $\tau(K)$:
\begin{equation}
\label{Birkhoff contraction coefficient}
\tau(K) := \sup_{0<h(\mu,\mu^\prime)<\infty} \frac{h(K\mu,K\mu^\prime)}{h(\mu,\mu^\prime)}.
\end{equation}
\end{definition}

Convergence of probability measures in the Hilbert projective metric is stronger than convergence in total variation. In fact, the following inequality holds (see \cite{Atar-Zeitouni}):
\begin{displaymath}
\|\mu - \mu^\prime\|_{tv} \leq \frac{2}{\log 3}h(\mu,\mu^\prime),\ \ {\rm for}\ \mu,\mu^\prime\in{\mathcal P}(E).
\end{displaymath}

The following result gives a uniform bound, in terms of the step errors, to the total variation distance of optimal filters corresponding to different parameters, one being equal to the true parameter value $\alpha$. The proof of this  estimate is given in \cite{Oudjane}. We have slightly altered the notation there, in order to make it consistent with the rest of this paper.

\begin{lemma}[\cite{Oudjane}, Cor. 4.5] 
\label{Oudjane1}
Suppose that the kernel $K_\alpha$ is mixing with some $\epsilon>0$. We define the step error with respect to the Hilbert projective metric as follows:
\begin{equation}
\label{h step error}
\delta^H_n(\theta,\alpha) = h(\Psi_n^\theta(\mu),{\mathbf K_n^\alpha}(\Psi_{n-1}^\theta(\mu))),
\end{equation}
where ${\mathbf K}_n^{\alpha}$ is the transition kernel at time $n$ of the optimal filter corresponding parameter $\alpha$, i.e.
\begin{equation}
\label{optimal kernel}
{\mathbf K}_n^{\alpha} (\mu)(dx) = \frac{\int_{E} g(Y_n - h(x)) K_{\alpha} (x^\prime,dx)\mu(dx^\prime)}{\int_{E}\int_{E} g(Y_n - h(z)) K_{\alpha} (z^\prime,dz)\mu(dz^\prime)}.
\end{equation}
Then, the total error is uniformly bounded in total variation by
\begin{equation}
\label{total variation error}
\sup_{n\geq 0}\| \Psi_n^\theta(\mu)-\Psi_n^\alpha(\mu) \|_{tv} \leq \frac{2}{\epsilon^2 \log 3}\sup_{n\geq 0}\delta^H_n(\theta,\alpha)
\end{equation}
where $\|\cdot\|_{tv}$ is the total variation norm on the space of measures.
\end{lemma}

Note that the quantities appearing in (\ref{total variation error}) are all random. It would have been more accurate to write 
\[ \sup_{n\geq 0}{\mathbb E}[\| \Psi_n^\theta(\mu)-\Psi_n^\alpha(\mu) \|_{tv}| Y_n,\dots,Y_1] \leq \frac{2}{\epsilon^2 \log 3}\sup_{n\geq 0}{\mathbb E}[\delta^H_n(\theta,\alpha)|Y_n,\dots,Y_1],\]
i.e. the inequality (\ref{total variation error}) holds for every fixed realization $(y_n,\dots,y_1)$ of the observation process. To simplify the notation, we will avoid writing this conditional expectation each time. The comparison between random quantities (depending on the observation process) that follow are meant to hold for each fixed realization.

The proof of lemma \ref{Oudjane1} is quite intuitive: the total error is written as the sum of local errors made when the kernels of the optimal filters are different only at one time step. Yet, the earlier the error occurs, the smaller the contribution of the local error to the total error, since the optimal filter ``corrects itself''. This self-correcting property is a consequence of the mixing property of the kernel, which guarantees that the contribution of the local errors to the total error is going to decay exponentially fast as time evolves and thus, the total error which is their sum is going to be bounded.

More specifically, the proof is based on decomposing the total error to the errors made at each time step and then using Birkhoff's contraction coefficient to get a bound for the total error with respect to the step errors (\ref{h step error}). That is, the error is decomposed to:
\begin{displaymath}
\Psi_n^\theta(\mu) - \Psi_n^\alpha(\mu) = \sum_{k=1}^n [{\mathbf K}^\alpha_{k+1,n}(\Psi_k^\theta(\mu)) - {\mathbf K}^\alpha_{k+1,n}({\mathbf K}^\alpha_k(\Psi_{k-1}^\theta(\mu)))],
\end{displaymath}
where ${\mathbf K}^\alpha_{k+1,n}={\mathbf K}^\alpha_{k+1}\circ\cdots\circ{\mathbf K}^\alpha_n$ is the $k$-to-$n$  transition kernel of the optimal filter corresponding to parameter $\alpha\in\Theta$.
Then, the following inequality that connects the total variation norm and the Hilbert projective metric is used:
\begin{displaymath}
\|K\mu - K\mu^\prime\|_{tv} \leq \frac{2}{\log 3}\tau(K) h(\mu,\mu^\prime),\ \ \ \forall K\in{\mathcal K}(E) {\rm\ and\ } \forall\mu\in{\mathcal P}(E),
\end{displaymath}
where $\tau$ is Birkhoff's contraction coefficient and ${\mathcal K}(E)$ is the space of transition kernels on $E$. The result follows by the fact that $\tau({\mathbf K}^\alpha_{k+1,n}) \leq (\frac{1-\epsilon^2}{1+\epsilon^2})^{n-k}$ (which is a consequence of the mixing property of $K_\alpha$ with mixing constant $\epsilon>0$). 

A similar result holds if we assume that the step errors are uniformly bounded with respect to the total variation norm instead of the Hilbert projective metric. Once again, the following estimate comes from \cite{Oudjane}. We rewrite it so that it fits in the setting of this paper.

\begin{lemma}[\ \cite{Oudjane}, Cor.4.7] 
\label{Oudjane2}
Suppose that the kernel $K_\alpha$ is mixing with some $\epsilon>0$. Let 
\begin{equation}
\label{tv step error}
\delta_n^{tv}(\theta,\alpha) := \|\Psi_n^\theta(\mu) - {\mathbf K_n^\alpha}(\Psi_{n-1}^\theta(\mu))\|_{tv},
\end{equation}
where ${\mathbf K_n^\alpha}$ is defined as in (\ref{optimal kernel}). Then, the total variation norm of the total error is uniformly bounded by
\begin{equation}
\label{tv norm error}
\sup_{n\geq 0}\|\Psi^\theta_n(\mu) - \Psi^\alpha_n(\mu)\|_{tv} \leq (1+\frac{2}{\epsilon^4\log3})\sup_{n\geq 0}\delta^{tv}_n(\theta,\alpha).
\end{equation}
\end{lemma}

The proof of the above lemma is based on the following estimate that provides an upper bound for the total variation norm in terms of the Hilbert projective norm: if the nonnegative kernel $K$ is mixing, then for any nonzero finite measures $\mu,\mu^\prime\in{\mathcal M}^+(E)$:
\[ h(K\mu,K\mu^\prime)\leq \frac{1}{\epsilon^2}\|\frac{\mu}{\mu(E)} - \frac{\mu^\prime}{\mu^\prime(E)}\|_{tv}. \]

It is now a straight-forward corollary that the uniform convergence of the step errors (\ref{h step error}) and (\ref{tv step error}) implies the uniform continuity of the optimal filters, in the total variation norm.

\begin{corollary}
\label{uniform continuity}
	Suppose that the Hilbert projective metric of the step errors converges to zero uniformly with respect to parameter $\theta$:
	\begin{equation}
	\label{h uniform convergence}
	\lim_{\theta\rightarrow\alpha} \sup_{n\geq 0} h(\Psi_n^\theta(\mu),{\mathbf K_n^\alpha}(\Psi_{n-1}^\theta(\mu))) = 0. 
	\end{equation}
	Then, the optimal filters are uniformly continuous with respect to $\theta$, in the total variation norm:
	\begin{equation}
	\label{tv uniform continuity}
	\lim_{\theta\rightarrow\alpha} \sup_{n\geq 0} \|\Psi_n^\theta(\mu) - \Psi_n^\alpha(\mu)\|_{tv} = 0.
	\end{equation}Uniform continuity of the optimal filter with respect to the parameter (\ref{tv uniform continuity}) also holds if 
	\begin{equation}
	\label{tv uniform convergence}
	\lim_{\theta\rightarrow\alpha}\sup_{n\geq 0}\|\Psi_n^\theta(\mu)-{\mathbf K_n^\alpha}(\Psi_{n-1}^\theta(\mu))\|_{tv}=0.
	\end{equation}
\end{corollary}

The problem now becomes how to show the uniform continuity of the step errors. It is easy to show the following

\begin{lemma}
Suppose that for every $\theta,\alpha\in\Theta$ and $x^\prime\in E$, the probability measures $K_\theta(x^\prime,\cdot)$ and $K_\alpha(x^\prime,\cdot)$ are absolutely continuous with respect to each other. If, in addition, there exist functions $c(\theta,\alpha)$ and $d(\theta,\alpha)$ such that 
\begin{equation}
\label{kernel uniform bound}
0 < c(\theta,\alpha)\leq \frac{dK_\theta(x^\prime,\cdot)}{dK_\alpha(x^\prime,\cdot)}(x)\leq c(\theta,\alpha)\exp(d(\theta,\alpha)) <+\infty
\end{equation}
then Hilbert projective distance of the step error is uniformly bounded by
\begin{equation*}
h(\Psi^\theta_n(\mu),{\mathbf K}_\alpha(\Psi^\theta_{n-1}(\mu))) \leq \log(\sup_{x^\prime,x} \frac{dK_\theta}{dK_\alpha}(x^\prime,x)\sup_{x^\prime,x} \frac{dK_\alpha}{dK_\theta}(x^\prime,x))=: h(K_\theta,K_\alpha) \leq d(\theta,\alpha).
\end{equation*}
Consequently, if $\lim_{\theta\rightarrow\alpha}d(\theta,\alpha)=0$, then the step error converges uniformly, i.e. (\ref{h uniform convergence}) holds.
\end{lemma}

The advantage of this lemma is that condition (\ref{kernel uniform bound}) only involves the kernels $K_\theta$ and $K_\alpha$ and can be easily checked, but it is too restrictive.  

The next lemma gives sufficient conditions for (\ref{tv uniform convergence}) to hold, that is less restrictive than (\ref{kernel uniform bound}), by making use of the Mean Value Theorem for real-valued functions defined on a Banach space. For the rest of this section, we assume that the parameter space $\Theta$ is a Banach space with norm $\|\cdot\|$.

\begin{lemma} 
\label{MVT}
Suppose that there exists a neighborhood ${\mathcal N}(\alpha)$ of $\alpha$, such that the Fr\'{e}chet derivative of kernel $K_\theta$ with respect to the weak operator norm exists for every $\theta\in{\mathcal N}(\alpha)$, i.e. $\exists\ L_\theta:{\mathcal C}_b(E)\rightarrow{\mathcal C}_b(E)$ such that
\begin{displaymath}
		\lim_{h\downarrow 0}|\frac{\mu (K_{\theta+h}f) - \mu (K_\theta f)}{\|h\|} - \mu (L_\theta f)| = 0,\ \forall \mu\in{\mathcal P}(E) \ {\rm and}\ f \in {\mathcal C}_b(E).
\end{displaymath}
Then, for $\theta$ sufficiently close to $\alpha$ ($\theta\in{\mathcal N}(\alpha)$),
\begin{equation}
\label{weak step bound}
|\Psi_n^\theta(\mu)(f) - {\mathbf K_n^\alpha}(\Psi_{n-1}^\theta(\mu))(f)| \leq 2\|f\|_{\infty}  \|\theta-\alpha\|\frac{\Psi_{n-1}^\theta(\mu)(|L_{\theta^\prime}g_n|)}{\Psi_{n-1}^\theta(\mu)(K_{\theta^\prime}g_n)},
\end{equation}
where by $g_n$ we denote the function $g_n(x) = g(Y_n-h(x))$. As an immediate consequence, we have that
\begin{equation}
\label{tv step bound}
\|\Psi_n^\theta(\mu) - {\mathbf K_n^\alpha}(\Psi_{n-1}^\theta(\mu))\|_{tv} \leq 2  \|\theta-\alpha\|\frac{\Psi_{n-1}^\theta(\mu)(|L_{\theta^\prime}g_n|)}{\Psi_{n-1}^\theta(\mu)(K_{\theta^\prime}g_n)}.
\end{equation}
\end{lemma}

\begin{proof}
For every $n>0$ and $\theta\in\Theta$, we define a real-valued function $F_{n,\theta}:\Theta\rightarrow{\mathbb R}$ by
\begin{displaymath}
F_{n,\theta}(\theta^\prime) = {\mathbf K}_n^{\theta^\prime}(\Psi_{n-1}^\theta(\mu))(f) = \frac{\Psi_{n-1}^\theta(\mu)(K_{\theta^\prime} f g_n)}{\Psi_{n-1}^\theta(\mu)(K_{\theta^\prime} g_n)},
\end{displaymath}
for some fixed $\theta\in{\mathcal N}(\alpha)$ and $f\in{\mathcal C}_b(E)$. Then, 
\begin{displaymath}
|\Psi_n^\theta(\mu)(f) - {\mathbf K_n^\alpha}(\Psi_{n-1}^\theta(\mu))(f)| =  |F_{n,\theta}(\theta) - F_{n,\theta}(\alpha)|.
\end{displaymath}

It is not hard to see that if the Fr\'{e}chet derivative of the kernel exists then the Fr\'{e}chet derivative $F^\prime_{n,\theta}$ of $F_{n,\theta}$ also exists and is equal to 
\begin{displaymath}
F_{n,\theta}^\prime(\theta^\prime) = \frac{\Psi_{n-1}^\theta(\mu)(L_{\theta^\prime} f g_n)\Psi_{n-1}^\theta(\mu)(K_{\theta^\prime} g_n) - \Psi_{n-1}^\theta(\mu)(K_{\theta^\prime} f g_n)\Psi_{n-1}^\theta(\mu)(L_{\theta^\prime} g_n)}{\Psi_{n-1}^\theta(\mu)(K_{\theta^\prime} g_n)^2},
\end{displaymath}
which can be easily bounded by:
\begin{eqnarray}
\label{derivative bound}
\nonumber &|F_{n,\theta}^\prime(\theta^\prime)| = \frac{|\int_{E^4} (f(x)-f(z))g_n(x)g_n(z)L_{\theta^\prime}(x^\prime,dx) K_{\theta^\prime}(z^\prime,dz) \Psi^\theta_{n-1}(\mu)(dx^\prime) \Psi^\theta_{n-1}(\mu)(dz^\prime)|} {\Psi_{n-1}^\theta(\mu)(K_{\theta^\prime} g_n)^2} \leq \\
\nonumber &\\
\nonumber &\leq \frac{\int_{E^3} \sup_x|f(x)-f(z)|\cdot|\int_E g_n(x)L_{\theta^\prime}(x^\prime,dx)|g_n(z) K_{\theta^\prime}(z^\prime,dz) \Psi^\theta_{n-1}(\mu)(dx^\prime) \Psi^\theta_{n-1}(\mu)(dz^\prime)} {\Psi_{n-1}^\theta(\mu)(K_{\theta^\prime} g_n)^2} \leq \\ \nonumber &\\
&\leq 2\|f\|_{\infty} \frac{\Psi_{n-1}^\theta(\mu)(|L_{\theta^\prime}g_n|)}{\Psi_{n-1}^\theta(\mu)(K_{\theta^\prime}g_n)}.
\end{eqnarray}
By the Mean Value Theorem (see, for example, \cite{Cheney}, p. 122), there exists a $\theta^\prime\in\{\alpha + t(\theta-\alpha): 0\leq t \leq 1\}$, such that
\begin{equation}
\label{MVT for F}
F_{n,\theta}(\theta) - F_{n,\theta}(\alpha) = F_{n,\theta}^\prime(\theta^\prime)(\theta-\alpha),
\end{equation}
Finally, (\ref{weak step bound}) and (\ref{tv step bound}) follow from (\ref{MVT for F}) and (\ref{derivative bound}).
\end{proof}

Usually, the parameter space $\Theta$ is Euclidean and, thus, the Fr\'{e}chet derivatives coincide with the ``usual derivatives'' of real-valued functions defined on the Euclidean space.
We can also get different bounds for the Fr\'{e}chet derivative. Suppose that the function 
\begin{equation}
\label{Lambda}
\Lambda_{\theta^\prime}(x) = \sup_{x^\prime} |\frac{dL_{\theta^\prime}(x^\prime,\cdot)}{dK_{\theta^\prime}(x^\prime,\cdot)}(x)| 
\end{equation}
is well defined. Then it is easy to see that
\begin{equation}
\label{derivative bound 2}
|F_{n,\theta}^\prime(\theta^\prime)| \leq 2\|f\|_\infty\frac{\Psi_{n-1}^\theta(\mu)(K_{\theta^\prime}(\Lambda_{\theta^\prime}g_n))}{\Psi_{n-1}^\theta(\mu)(K_{\theta^\prime}g_n)} = 2\|f\|_\infty{\mathbb E}_{\mu,(n,\theta,\theta^\prime)}(\Lambda_{\theta^\prime}(X_n)|Y_n,\dots,Y_1),
\end{equation}
where we denote by ${\mathbb P}_{\mu,(n,\theta,\theta^\prime)}$ the distribution of a Markov chain with initial distribution $\mu$ and transition kernel $K_\theta$ up to time $n-1$ and $K_{\theta^\prime}$ at times $\geq n$ and by ${\mathbb E}_{\mu,(n,\theta,\theta^\prime)}$ the respective expectation. 

So, the uniform boundedness of (\ref{derivative bound}) or (\ref{derivative bound 2}) would imply the uniform continuity of the optimal filters with respect to the parameter, in the total variation norm. Working with (\ref{derivative bound 2}) can be more intuitive, even though it is actually stronger than (\ref{derivative bound}). Note, though, that the function $\Lambda_{\theta^\prime}$ is not assumed to be bounded (this would have been too restrictive), so getting a uniform bound for (\ref{derivative bound 2}), for any given observation process, is not always possible. However, if we assume that the state space $E$ is compact then the continuity of $|\frac{dL_{\theta^\prime}(x^\prime,\cdot)}{dK_{\theta^\prime}(x^\prime,\cdot)}(x)|$ with respect to $x$ and $x^\prime$ would imply the uniform continuity of the optimal filters, in the sense of (\ref{tv uniform continuity}).

The above discussion is summarized in the following

\begin{corollary}
\label{conditions for uniform continuity}
Under the mixing assumption for $K_\alpha$ and for $\theta$ sufficiently close to $\alpha$ ($\theta\in{\mathcal N}(\alpha)$), there exists a $\theta^\prime\in{\mathcal N}(\alpha)$ so that
\begin{displaymath}
\|\Psi^\theta_n(\mu)-\Psi^\alpha_n(\mu)\|_{tv}\leq 2(1+\frac{2}{\epsilon^4 \log3})\|\theta-\alpha\|\frac{\Psi_{n-1}^\theta(\mu)(|L_{\theta^\prime}g_n|)}{\Psi_{n-1}^\theta(\mu)(K_{\theta^\prime}g_n)}.
\end{displaymath}
If the function $\Lambda_{\theta^\prime}$ is well defined, then 
\begin{displaymath}
\|\Psi^\theta_n(\mu)-\Psi^\alpha_n(\mu)\|_{tv}\leq 2(1+\frac{2}{\epsilon^4 \log3})\|\theta-\alpha\|{\mathbb E}_{\mu,(n,\theta,\theta^\prime)}(\Lambda_{\theta^\prime}(X_n)|Y_n,\dots,Y_1).
\end{displaymath}
Consequently, if there exists an $M$ such that 
\begin{equation}
\label{bounded conditional mean}
\sup_{n>0}\sup_{\theta,\theta^\prime\in{\mathcal N}(\alpha)} \frac{\Psi_{n-1}^\theta(\mu)(|L_{\theta^\prime}g_n|)}{\Psi_{n-1}^\theta(\mu)(K_{\theta^\prime}g_n)}(Y_1(\omega),\dots,Y_n(\omega))<M
\end{equation}
for any realization $\omega\in\Omega$ and a sufficiently small neighborhood ${\mathcal N}(\alpha)$ of $\alpha$, then (\ref{tv uniform continuity}) holds.
\end{corollary}

Note, however, that for the asymptotic stability of the optimal filter to hold it is sufficient to show the almost sure uniform continuity of the optimal filters, i.e.
\begin{equation}
\label{a.s. tv uniform continuity}
\lim_{\theta\rightarrow\alpha}\sup_{n>0}{\mathbb E}_{\mu^\prime,\alpha}\|\Psi^n_\theta(\mu)-\Psi^n_\alpha(\mu)\|_{tv} = 0.
\end{equation}
As before, for this to be true it is sufficient to show that
\[ \sup_{n>0}\sup_{\theta,\theta^\prime\in{\mathcal N}(\alpha)} {\mathbb E}_{\mu^\prime,\alpha}\frac{\Psi_{n-1}^\theta(\mu)(|L_{\theta^\prime}g_n|)}{\Psi_{n-1}^\theta(\mu)(K_{\theta^\prime}g_n)} <+\infty \]
or, if $\Lambda_{\theta^\prime}$ is well defined, that
\begin{equation}
\sup_{n>0}\sup_{\theta,\theta^\prime\in{\mathcal N}(\alpha)} {\mathbb E}_{\mu^\prime,\alpha}{\mathbb E}_{\mu,(n,\theta,\theta^\prime)}(\Lambda_{\theta^\prime}(X_n)|Y_n,\dots,Y_1) < \infty.
\end{equation}
Again, this last condition is much easier to work with. We write
\begin{eqnarray}
\nonumber &{\mathbb E}_{\mu^\prime,\alpha}{\mathbb E}_{\mu,(n,\theta,\theta^\prime)}(\Lambda_{\theta^\prime}(X_n)|Y_n,\dots,Y_1) = {\mathbb E}_{\mu,(n,\theta,\theta^\prime)}(\Lambda_{\theta^\prime}(X_n)\cdot\frac{dQ^n_{\mu^\prime,\alpha}} {dQ^n_{\mu,(n,\theta,\theta^\prime)}}(Y_n,\dots,Y_1)) \leq \\
\nonumber &\leq ({\mathbb E}_{\mu,(n,\theta,\theta^\prime)}\Lambda_{\theta^\prime}(X_{n-1})^{n+1})^{\frac{1}{n+1}}\cdot ({\mathbb E}_{\mu^\prime,\alpha}(\frac{dQ^n_{\mu^\prime,\alpha}} {dQ^n_{\mu,(n,\theta,\theta^\prime)}}(Y_n,\dots,Y_1))^{\frac{1}{n}})^{\frac{n}{n+1}} = \\
&= ({\mathbb E}_{\mu,\theta}K_{\theta^\prime}\Lambda_{\theta^\prime}(X_{n-1})^{n+1})^{\frac{1}{n+1}}\cdot ({\mathbb E}_{\mu^\prime,\alpha}(\frac{dQ^n_{\mu^\prime,\alpha}} {dQ^n_{\mu,(n,\theta,\theta^\prime)}}(Y_n,\dots,Y_1))^{\frac{1}{n}})^{\frac{n}{n+1}}.
\end{eqnarray}
It is easy to see that 
\begin{eqnarray*}
&\frac{dQ^n_{\mu^\prime,\alpha}} {dQ^n_{\mu,(n,\theta,\theta^\prime)}}(Y_n,\dots,Y_1)) = \frac{{\mathbb E}_{\mu^\prime,\alpha}(e^{-\frac{1}{2}\sum_{k=1}^n(Y_k-h(X_k))^2}|Y_n,\dots,Y_1)}{{\mathbb E}_{\mu,(n,\theta,\theta^\prime)}(e^{-\frac{1}{2}\sum_{k=1}^n(Y_k-h(X_k))^2}|Y_n,\dots,Y_1)}\leq \\
& \leq e^{\frac{1}{2}\sum_{k=1}^n{\mathbb E}_{\mu,(n,\theta,\theta^\prime)}((Y_k-h(X_k))^2|Y_n,\dots,Y_1)},
\end{eqnarray*}
by Jensen's inequality and the fact that $e^{-\frac{1}{2}\sum_{k=1}^n(Y_k-h(X_k))^2}\leq 1$. Consequently,
\begin{equation*}
{\mathbb E}_{\mu^\prime,\alpha}(\frac{dQ^n_{\mu^\prime,\alpha}} {dQ^n_{\mu,(n,\theta,\theta^\prime)}}(Y_n,\dots,Y_1))^{\frac{1}{n}} \leq {\mathbb E}_{\mu^\prime,\alpha}e^{\frac{1}{2n}\sum_{k=1}^n{\mathbb E}_{\mu,(n,\theta,\theta^\prime)}((Y_k-h(X_k))^2|Y_n,\dots,Y_1)},
\end{equation*}
But since the observation process is ergodic and function $h$ is bounded, this is also going to be bounded, provided that the two first moments of the limiting distribution $\nu_\alpha$ exist. This proves the following

\begin{lemma}
\label{a.s. uniform continuity lemma}
Suppose that the function $\Lambda_\theta(x)$ defined in (\ref{Lambda}) is well-defined and 
\begin{equation}
\label{simple condition}
\sup_{n>0}\sup_{\theta,\theta^\prime\in{\mathcal N}(\alpha)}({\mathbb E}_{\mu,\theta}K_{\theta^\prime}\Lambda_{\theta^\prime}(X_{n-1})^{n+1})^{\frac{1}{n+1}} <+\infty,
\end{equation}
for some sufficiently small neighborhood of $\alpha$, ${\mathcal N}(\alpha)$. Then, if $K_\alpha$ is mixing and the two first moments of $\nu_\alpha$ exist, the optimal filters will almost surely be uniformly continuous with respect to the parameters, i.e. (\ref{a.s. tv uniform continuity}) will hold.
\end{lemma}

The importance of this lemma is that its assumptions depend only on properties of the kernels. We can now prove the asymptotic stability of the optimal filter of system \ref{system B}, with respect to the initial conditions.

\section{Asymptotic stability of non-mixing systems}

In general, the asymptotic behavior of the optimal filter is not well understood, even when the system is ergodic (see \cite{Lipster} and \cite{Budhiraja}, regarding the gap in the proof of the ergodicity of the optimal filters, in \cite{Kunita 71}, and the consequences). If the state space is finite (\cite{Lipster}), or the kernel is mixing (\cite{Atar-Zeitouni}), it has been shown that the optimal filters will be asymptotically stable. In fact, Chigansky and Lipster, in \cite{Chigansky}, recently showed that the optimal filter will also be stable under conditions that are weaker than mixing but stronger than ergodicity of the signal.

In the following lemma, we prove the asymptotic stability of the optimal filter for System \ref{system B}, which is a non-mixing system. 

\begin{lemma} 
\label{main lemma}
Let $Y$ be as in system \ref{system A} and suppose that the assumptions of lemma \ref{posterior lemma 1} are satisfied. We also assume the uniform continuity of the optimal filters, in the sense of (\ref{a.s. tv uniform continuity}). Then 
\begin{equation}
\label{tv marginal stability}
  \lim_{n\rightarrow\infty}{\mathbb E}_{\mu^\prime,\alpha}\|\Psi^u_n(\mu)-\Psi^\alpha_n(\mu^\prime)\|_{tv}=0,
\end{equation}
for any initial distributions $\mu,\mu^\prime \in {\mathcal P}(E)$.
\end{lemma}

\begin{proof}
First, we note that it is sufficient to show (\ref{tv marginal stability}) for the same initial distributions, since the rest follows by the asymptotic stability of the optimal filters, for mixing kernels. That is, we want to show 
\begin{equation*}
  \lim_{n\rightarrow\infty}{\mathbb E}_{\mu^\prime,\alpha}\|\Psi^u_n(\mu)-\Psi^\alpha_n(\mu)\|_{tv}=0,
\end{equation*}
We decompose the optimal filters as follows:
\begin{eqnarray*}
\Phi_n(\mu\otimes u)(dx,d\theta) &=& \tilde{\mathbb P}_{\mu,u}(X_n\in dx,\theta_n\in d\theta|Y_n,\ldots,Y_1) = \\
&=& \tilde{\mathbb P}_{\mu,u}(X_n\in dx | Y_n,\ldots,Y_1,\theta_n=\theta)\tilde{\mathbb P}_{\mu,u}(\theta_n\in d\theta|Y_n,\ldots,Y_1) = \\
&=& {\mathbb P}_{\mu,\theta}(X_n\in dx|Y_n,\ldots,Y_1)\tilde{\mathbb P}_{\mu,u}(d\theta|Y_n,\ldots,Y_1) = \\
&=& \Psi_n^\theta(\mu)(dx)Z^{\mu,u}_n(d\theta),
\end{eqnarray*}
where $Z_n^{\mu,u}(d\theta)=\tilde{\mathbb P}_{\mu,u}(\theta\in d\theta|Y_n,\ldots,Y_1)$. So,
\begin{eqnarray}
\nonumber {\mathbb E}_{\mu^\prime,\alpha}\|\Psi^u_n(\mu) - \Psi^\alpha_n(\mu)\|_{tv} &=& {\mathbb E}_{\mu^\prime,\alpha} \|\int_\Theta\Psi_n^\theta(\mu)Z^{\mu,u}_n(d\theta) - \Psi_n^{\alpha}(\mu)\|_{tv} = \\
\nonumber &=& {\mathbb E}_{\mu^\prime,\alpha}\|\int_\Theta(\Psi_n^\theta(\mu) - \Psi_n^{\alpha}(\mu))Z^{\mu,u}_n(d\theta)\|_{tv} \leq \\  
& \leq & {\mathbb E}_{\mu^\prime,\alpha}\int_\Theta \|\Psi_n^\theta(\mu) - \Psi_n^{\alpha}(\mu)\|_{tv} Z^{\mu,u}_n(d\theta).
\end{eqnarray}
Since we have assumed the uniform continuity of the optimal filters (\ref{a.s. tv uniform continuity}), $\forall\epsilon>0$ we can find a neighborhood ${\mathcal N}_\eta(\alpha)$ of $\alpha$ for some $\eta>0$, such that 
\begin{displaymath}
\forall n\geq 0,\ \ \ \sup_{\theta\in{\mathcal N}_\eta(\alpha)}{\mathbb E}_{\mu^\prime,\alpha}\|\Psi_n^\theta(\mu) - \Psi_n^{\alpha}(\mu)\|_{tv}<\frac{\epsilon}{2}.
\end{displaymath}
Also, by lemma \ref{posterior lemma 1}, we can find an $n_0$ such that
\begin{displaymath}
\forall n\geq n_0,\ \ {\mathbb E}_{\mu^\prime,\alpha}Z^{\mu,u}_n({\mathcal N}_\eta(\alpha)^c) <\frac{\epsilon}{4},
\end{displaymath} 
where ${\mathcal N}_\eta(\alpha)^c$ is the complement of ${\mathcal N}_\eta(\alpha)$. So, putting the last two estimates together, we get that $\forall n\geq n_0$
\begin{eqnarray*}
&{\mathbb E}_{\mu^\prime,\alpha}\|\Psi^u_n(\mu) - \Psi^\alpha_n(\mu)\|_{tv} \leq \\
&\sup_{\theta\in{\mathcal N}_\eta(\alpha)}{\mathbb E}_{\mu^\prime,\alpha}\|\Psi_n^\theta(\mu) - \Psi_n^{\alpha}(\mu)\|_{tv} + 2{\mathbb E}_{\mu^\prime,\alpha}Z^{\mu,u}_n({\mathcal N}_\eta(\alpha)^c) <\epsilon,
\end{eqnarray*}
which proves (\ref{tv marginal stability}).
\end{proof}

In the following theorem, we prove the asymptotic stability of the optimal filters of system \ref{system B}, under assumption that only involve the kernels $K_\theta$;

\begin{theorem}
Suppose that $Y$ is as in system \ref{system A} and that the parameter space is a compact Banach space. We also assume the following
\begin{enumerate}
	\item{ the prior distribution $u$ on the parameter space is such that there exist a sequence $\epsilon_n\downarrow0$ and a function $p:{\mathbb N}\rightarrow[1,\infty)$ satisfying (\ref{p(n)}), so that 
\begin{equation*}
\lim_{n\uparrow\infty}{\mathbb E}_{\mu,\alpha}[(\frac{\sup_{\theta\in{\mathcal N}_{\epsilon_n}(\alpha)}\frac{dQ^n_{\mu,\alpha}}{dQ^n_{\mu^\prime,\theta}}(Y_n,\dots,Y_1)}{u({\mathcal N}_{\epsilon_n}(\alpha))})^{\frac{1}{p(n)}}]<+\infty.
\end{equation*}}
	\item{for each $\theta\in\Theta$, the observation process $\{Y_n\}_{n>0}$ satisfies the LDP with rate function $J_\theta$. }
	\item{The mapping $\Theta\ni\theta\mapsto{\nu_\theta}\in{\mathcal P}({\mathbb R}^p)$, is an open mapping, in the sense of (\ref{open mapping}).}
	\item{The kernel $K_\theta$ and the rate function $J_\theta$ are continuous with respect to $\theta$.}
	\item{The function $\Lambda_\theta(x)$ defined in (\ref{Lambda}) is well-defined and 
\begin{equation*}
\sup_{n>0}\sup_{\theta,\theta^\prime\in{\mathcal N}(\alpha)}({\mathbb E}_{\mu,\theta}K_{\theta^\prime}\Lambda_{\theta^\prime}(X_{n-1})^{n+1})^{\frac{1}{n+1}} <+\infty,
\end{equation*}
for some sufficiently small neighborhood of $\alpha$, ${\mathcal N}(\alpha)$.}
	\item{The first two moments of $\nu_\alpha$ exist.}
\end{enumerate}
Then, the optimal filter of system \ref{system B} will eventually correct itself, i.e. it will satisfy (\ref{tv marginal stability}).
\end{theorem}

\begin{proof}
Just combine lemmas \ref{posterior lemma 2}, \ref{a.s. uniform continuity lemma} and \ref{main lemma}.
\end{proof}



\section*{Acknowledgments}
I would like to thank my PhD adviser, Ren\'{e} Carmona, for his help and support. It is his insights and his questions that inspired this work. My gratitude also goes to Pierre Del Moral and Ioannis Karatzas, for their useful comments and suggestions after reading earlier versions of this paper and to the anonymous referee for his/her invaluable advice. Finally, I would like to thank Yannis Kevrekidis who has been my host during these last months I have been working on this paper.

\end{document}